\documentclass[12pt]{amsart}

\usepackage{times,epsf,amssymb,amsmath,hyperref, wrapfig, rlepsf, color}

\begin{document}

\newtheorem{thm}{Theorem}[section]
\newtheorem{lem}[thm]{Lemma}
\newtheorem{cor}[thm]{Corollary}
\newtheorem{conj}[thm]{Conjecture}

\theoremstyle{definition}
\newtheorem{defn}[thm]{\bf{Definition}}

\theoremstyle{remark}
\newtheorem{rmk}[thm]{Remark}
\newtheorem{question}[thm]{Question}

\def\square{\hfill${\vcenter{\vbox{\hrule height.4pt \hbox{\vrule width.4pt height7pt \kern7pt \vrule width.4pt} \hrule height.4pt}}}$}

\newcommand{\SI}{\partial_\infty ({\Bbb H}^2\times {\Bbb R})}
\newcommand{\Si}{S^1_{\infty}\times {\Bbb R}}
\newcommand{\si}{S^1_{\infty}}
\newcommand{\PI}{\partial_{\infty}}

\newcommand{\BH}{\Bbb H}
\newcommand{\BHH}{{\Bbb H}^2\times {\Bbb R}}
\newcommand{\BR}{\Bbb R}
\newcommand{\BC}{\Bbb C}
\newcommand{\BZ}{\Bbb Z}
\newcommand{\caps}{\overline{\BH^2}\times\{\pm\infty\}}

\newcommand{\e}{\epsilon}

\newcommand{\wh}{\widehat}
\newcommand{\wt}{\widetilde}

\newcommand{\A}{\mathcal{A}}
\newcommand{\C}{\mathcal{C}}
\newcommand{\p}{\mathcal{P}}
\newcommand{\R}{\mathcal{R}}
\newcommand{\B}{\mathcal{B}}
\newcommand{\h}{\mathcal{H}}
\newcommand{\T}{\mathfrak{T}}
\newcommand{\s}{\mathcal{S}}
\newcommand{\U}{\mathcal{U}}
\newcommand{\V}{\mathcal{V}}
\newcommand{\X}{\mathcal{X}}
\newcommand{\tc}{\mathcal{T}}
\newcommand{\Y}{\mathbf{Y}}

\newenvironment{pf}{{\it Proof:}\quad}{\square \vskip 12pt}

\title[Asymptotic Plateau Problem in $\BHH$]{Asymptotic Plateau Problem in $\BHH$: Tall Curves}
\author{Baris Coskunuzer}
\address{UT Dallas, Dept. Math. Sciences, Richardson, TX 75080}
\email{coskunuz@utdallas.edu}
\thanks{The author is partially supported by Simons Collaboration Grant, and Royal Society Newton Mobility Grant.}

\maketitle

\begin{abstract} We study the asymptotic Plateau problem in $\BHH$ for area minimizing surfaces, and give a fairly complete solution for finite curves. 
\end{abstract}

\section{Introduction}

Asymptotic Plateau Problem in $\BHH$ asks the existence of a minimal surface $\Sigma$ in $\BHH$ for a given curve $\Gamma$ in $\SI$ with $\PI\Sigma=\Gamma$. 
In the past years, the existence, uniqueness, and regularity of solutions to the asymptotic Plateau problem in $\BHH$ have been studied extensively, and many important results have been  obtained by the leading researchers of the field, e.g. \cite{ CR, CMT, Da, FMMR,  KM, MMR, MoR, MRR, NR, PR, RT, ST, ST2}.


Unlike $\BH^3$, the asymptotic Plateau problem in $\BHH$ is quite interesting and challenging as  there are several families of curves in $\Si$, which does not bound any minimal surface in $\BHH$ \cite{ST}. In this paper, we finish off an important case  by classifying strongly fillable, finite curves in $\Si$ as follows:

\begin{thm} \label{main2} Let $\Gamma$ be a finite collection of disjoint, smooth Jordan curves in $\Si$ with $h(\Gamma)\neq \pi$. Then, there exists an area minimizing surface $\Sigma$ in $\BHH$ with $\PI \Sigma = \Gamma$ if and only if $\Gamma$ is a tall curve.
\end{thm}


The organization of the paper is as follows. In the next section, we give some definitions, and introduce the basic tools which we use in our construction. In Section \ref{tallcurvesec}, we introduce tall curves, and study their properties. In Section \ref{secAPP}, we prove our main result above. In Section \ref{secAPPmin}, we show that the asymptotic Plateau problem for minimal surfaces and area minimizing surfaces are quite different, and construct some explicit examples. Finally in Section \ref{secremarks}, we give some concluding remarks, and mention some interesting open problems in the subject. We postpone some technical steps to the appendix at the end.

\subsection{Acknowledgements}
Part of this research was carried out at MIT and Max-Planck Institute during our visit. The author would like to thank them for their great hospitality.
An earlier version of this paper appeared as \cite{Co1} 

 

\section{Preliminaries} \label{secprem}

In this section, we give the basic definitions, and a brief overview of the past results which will be used in the paper. 

Throughout the paper, we use the product compactification of $\BHH$. In particular, $\overline{\BHH}=\overline{\BH^2}\times\overline{\BR}=\BHH\cup \SI$ where $\SI$ consists of three components, i.e. the infinite open cylinder $\Si$ and the closed caps at infinity $\overline{\BH^2}\times\{+\infty\}$, $\overline{\BH^2}\times\{-\infty\}$. Hence, $\overline{\BHH}$ is a solid cylinder under this compactification.

Let $\Sigma$ be an open, complete surface in $\BHH$, and $\PI\Sigma$ represent the asymptotic boundary of $\Sigma$ in $\SI$. Then, if $\overline{\Sigma}$ is the closure of $\Sigma$ in $\overline{\BHH}$, then $\PI \Sigma= \overline{\Sigma}\cap \SI$.

\begin{defn} A surface is {\em minimal} if the mean curvature $H$ vanishes everywhere. A compact surface with boundary $\Sigma$ is called {\em area minimizing surface} if $\Sigma$ has the smallest area among the surfaces with the same boundary. A noncompact surface is called {\em area minimizing surface} if any compact subsurface is an area minimizing surface. 
\end{defn}

In this paper, we study the Jordan curves in $\SI$ which bounds a complete, embedded, minimal surfaces in $\BHH$.
Throughout the paper, when we say {\em a curve in $\SI$} we mean a finite collection of pairwise disjoint Jordan curves in $\SI$.

\begin{defn} [Fillable Curves] Let $\Gamma$ be a curve in $\SI$. We call $\Gamma$ {\em fillable} if $\Gamma$ bounds a complete, embedded, minimal surface $S$ in $\BHH$, i.e. $\PI S=\Gamma$. We call $\Gamma$ {\em strongly fillable} if $\Gamma$ bounds a complete, embedded, area minimizing surface $\Sigma$ in $\BHH$, i.e. $\PI \Sigma=\Gamma$. 
\end{defn}

Notice that a strongly fillable curve is fillable since any area minimizing surface is minimal. 

\begin{defn} [Finite and Infinite Curves] Let $\Gamma$ be a curve in $\SI$.  Decompose $\Gamma=\Gamma^+\cup\Gamma^-\cup\wt{\Gamma}$ such that 	$\Gamma^\pm=\Gamma\cap (\overline{\BH^2}\times\{\pm\infty\})$ and $\wt{\Gamma}=\Gamma\cap(\Si)$. In particular, $\Gamma^\pm$ is a collection of closed arcs and points in the closed caps at infinity, where $\wt{\Gamma}$ is a collection of open arcs and closed curves in the infinite open cylinder. With this notation, we call a curve $\Gamma$ {\em finite} if $\Gamma^+=\Gamma^-=\emptyset$. We call $\Gamma$ {\em infinite} otherwise.
\end{defn}

\noindent \textbf{Asymptotic Plateau Problem for $\BHH$:} 

\vspace{.2cm}

\noindent \textit{Which (finite or infinite) $\Gamma$ in $\SI$ is fillable or strongly fillable?}

As the question suggests, there are mainly four versions of the problem: Classification of "Fillable finite curves", "Fillable infinite curves",  "Strongly fillable finite curves", and "Strongly fillable infinite curves". Unfortunately, we are currently far from classification of the fillable (finite or infinite) curves \cite{FMMR}.

Recently, we gave a classification for strongly fillable, infinite curves in \cite{Co2}. In this paper, we give a fairly complete solution for the classification of {\em strongly fillable, finite curves} in $\SI$.

One of the most interesting properties of the asymptotic Plateau problem in $\BHH$ is the existence of non-fillable curves. While any curve $\Lambda$ in $S^2_\infty(\BH^3)$ is strongly fillable in $\BH^3$ \cite{An}, Sa Earp and Toubiana showed that there are some non-fillable $\Gamma$ in $\SI$ \cite{ST}.

\begin{defn} \label{taildef} [Thin tail] Let $\Gamma$ be a Jordan curve in $\SI$, and let $\tau$ be an arc in $\Gamma$. Assume that there is a vertical straight line $L_0$ in $\Si$ such that
	
	\begin{itemize}
		
		\item $\tau \cap L_0 \neq \emptyset$ and $\partial \tau\cap L_0 = \emptyset$,
		
		\item $\tau$ stays in one side of $L_0$,
		
		\item $\tau\subset \si \times (c,c+\pi)$ for some $c\in \BR$.
		
	\end{itemize}
	
	Then, we call $\tau$ {\em a thin tail} in $\Gamma$.
	
\end{defn}

\begin{lem} \cite{ST} \label{thinlem} Let $\Gamma$ be a curve in $\SI$. If $\Gamma$ contains a thin tail, then there is no properly immersed minimal surface $\Sigma$ in $\BHH$ with $\PI \Sigma=\Gamma$.
\end{lem}

When $\Gamma$ is an essential smooth simple closed curve in $\Si$ which is a vertical graph over $\si\times\{0\}$, then the vertical graphs over $\BH^2$ gives a positive answer to this existence question \cite{NR}. However, for nonessential (nullhomotopic) simple closed curves in $\Si$, the situation is quite different. Unlike the $\BH^3$ case \cite{An}, Sa Earp and Toubiana showed that there are some simple closed curves $\Gamma$ in $\Si$ which are not fillable \cite{ST}.

\begin{defn} \label{tail} [Thin Tail] Let $\gamma$ be an arc in $\Si$. Assume that there is a vertical straight line $L_0$ in $\Si$ such that

\begin{itemize}

\item $\gamma \cap L_0 \neq \emptyset$ and $\partial \gamma\cap L_0 = \emptyset$,

\item $\gamma$ stays in one side of $L_0$,

\item $\gamma\subset \si \times (c,c+\pi)$ for some $c\in \BR$.

\end{itemize}

Then, we call $\gamma$ {\em a thin tail} in $\Gamma$.

\end{defn}

With the definition above, we have the following nonexistence result: 

\begin{lem} \cite{ST} \label{nonexist} Let $\Gamma$ be a simple closed curve in $\SI$. If $\Gamma$ contains a thin tail, then there is no properly immersed minimal surface $\Sigma$ in $\BHH$ with $\PI \Sigma=\Gamma$.
\end{lem}

This nonexistence result makes the asymptotic Plateau problem quite interesting. In particular, to address the fillability question, we need to understand which curves have no thin tails. In Section \ref{tallcurvesec}, we will introduce a notion called {\em tall curves} to recognize them.

To construct our sequence of compact area minimizing surfaces in our main result, we need the following classical result of geometric measure theory.

\begin{lem} \label{AMSexist} \cite[Theorem 5.1.6 and 5.4.7]{Fe} [Existence and Regularity of Area Minimizing Surfaces] Let $M$ be a homogeneously regular, closed (or mean convex) $3$-manifold. Let $\gamma$ be a nullhomologous smooth curve in $\gamma$. Then, $\gamma$ bounds an area minimizing surface $\Sigma$ in $M$. Furthermore, any such area minimizing surface is smoothly embedded.	
\end{lem}

Now, we state the convergence theorem for area minimizing surfaces, which will be used throughout the paper. Note that we use convergence in the sense of Geometric Measure Theory, i.e. the convergence of rectifiable currents in the flat metric.

\begin{lem} \label{convergence} [Convergence] Let $\{\Sigma_i\}$ be a sequence of complete area minimizing surfaces in $\BHH$ where $\Gamma_i= \PI \Sigma_i$ is a finite collection of closed curves in $\Si$. If $\Gamma_i$ converges to a finite collection of closed curves $\wh{\Gamma}$ in $\Si$, then there exists a subsequence $\{\Sigma_{n_j}\}$ such that $\Sigma_{n_j}$ converges to an area minimizing surface $\wh{\Sigma}$ (possibly empty) with $\PI\wh{\Sigma}\subset\wh{\Gamma}$. In particular, the convergence is smooth on compact subsets of $\BHH$.
\end{lem}

\begin{pf}  Let $\Delta_n = \mathbf{B}_n(0)\times [-C,C]$ be convex domains in $\BHH$ where $\mathbf{B}_n(0)$ is the closed disk of radius $n$ in $\BH^2$ with center $0$, and $\wh{\Gamma}\subset \si\times (-C,C)$. For $n$ sufficiently large, consider the surfaces $S_i^n=\Sigma_i\cap \Delta_n$. Since the area of the surfaces $\{S^n_i\subset \Delta_n\}$ is uniformly bounded by $|\partial \Delta_n|$, and $\partial S^n_i$ can be bounded by using standard techniques. Hence, if $\{S^n_i\}$ is an infinite sequence, then we get a convergent subsequence of $\{S^n_i\}$ in $\Delta_n$ with {\em nonempty limit} $S^n$. $S^n$ is an area minimizing surface in $\Delta_n$ by the compactness theorem  for rectifiable currents (codimension-1) with the flat metric of Geometric Measure Theory \cite{Fe}. By the regularity theory, the limit $S^n$ is a smoothly embedded area minimizing surface in $\Delta^n$.
	
If the sequence $\{S^n_i\}$ is an infinite sequence for infinitely many $n$, we get an infinite sequence of compact area minimizing surfaces $\{S^n\}$. Then, by using the diagonal sequence argument, we can find a subsequence of $\{\Sigma_i\}$ converging to an area minimizing surface $\wh{\Sigma}$ with $\PI \wh{\Sigma}\subset \wh{\Gamma}$ as $\Gamma_i\to \wh{\Gamma}$. Note also that for fixed $n$, the curvatures of $\{S^n_i\}$ are uniformly bounded by curvature estimates for area minimizing surfaces. Hence, with the uniform area bound, we get smooth convergence on compact subsets of $\BHH$. For further details, see \cite[Theorem 3.3]{MW}.
\end{pf}

\begin{rmk} \label{escape2} (Empty Limit) In the proof above, there might be cases like $\{S^n_i\}$ is a finite sequence for any $n$. In particular, assume that if for every $n$, there exists $K_n>>0$ such that for every $i>K_n$, $\Sigma_i\cap \Delta_n=\emptyset$. In such a case the limit is empty, and we say $\{\Sigma_i\}$ {\em escapes to infinity}. An example to this case is a sequence of rectangles $R_i$ in $\Si$ with $h(R_i)\searrow \pi$ and $R_i\to\wh{R}$ where $\wh{R}$ is a rectangle of height $\pi$. Then, the sequence of area minimizing surfaces $P_i$ with $\PI P_i=R_i$ escapes to infinity, as there is no area minimizing surface $\Sigma$ with $\PI \Sigma = \wh{\Gamma}$. In Theorem \ref{APP}, we will prove that if $\wh{\Gamma}$ is a tall curve, the sequence $\{\Sigma_i\}$ does not escape to infinity, and a subsequence $\Sigma_{i_j}$ converges to an area minimizing surface $\wh{\Sigma}$ with $\PI\wh{\Sigma}\subset\wh{\Gamma}$.
\end{rmk}

\begin{rmk} [Asymptotic Regularity] In Lemma \ref{surface}, we proved $\C^0$ asymptotic regularity for area minimizing surfaces bounding tall curves. Note that Kloneckner and Mazzeo proved higher order asymptotic regularity for embedded minimal surfaces in $\BHH$ \cite[Section 3]{KM}.
\end{rmk}

\section{Tall Curves in $\Si$} \label{tallcurvesec} 

After Sa Earp - Toubiana's nonexistence result (Lemma \ref{nonexist}), one needs to understand the curves with no thin tails in order to solve asymptotic Plateau problem. In this section, we introduce a notion called {\em tall curves} to easily identify such curves. First, we study the tall rectangles. Then, by using these, we define the tall curves.

\vspace{.2cm}

\subsection{Tall Rectangles} \label{tallrecsec}

\vspace{.2cm}

\begin{defn} \label{tallrecdef} [Tall Rectangles] Consider the asymptotic cylinder $\Si$ with the coordinates $(\theta, t)$ where $\theta\in [0,2\pi)$ and $t\in \BR$. We call a rectangle $R=[\theta_1,\theta_2]\times[t_1,t_2]\subset \Si$  {\em tall rectangle} if $t_2-t_1>\pi$.
\end{defn}

In \cite{ST}, for the boundaries of tall rectangles, Sa Earp and Toubiana further proved the following:

\begin{lem} \cite{ST} \label{stdisk} If $R$ is a tall rectangle in $\Si$, then there exists a minimal surface $P$ in $\BHH$ with $\PI P = \partial R$. In particular, $P$ is a graph over $R$.
\end{lem}

Furthermore, they gave a very explicit description of $P$ as follows. Without loss of generality, let $R=[-\theta_1,\theta_1]\times[-c,c]$ in $\Si$ where $c>\frac{\pi}{2}$ and $\theta_1\in (0,\pi)$. Let $\varphi_t$ be the hyperbolic isometry of $\BH^2$ fixing the geodesic $\gamma$ with $\PI\gamma=\{-\theta_1,\theta_1\}$ with translation length $t$. Let $\wh{\varphi}$ be the isometry of $\BHH$ with $\wh{\varphi}_t(q,z)=(\varphi_t(q), z)$. They proved that $P$ is invariant under $\wh{\varphi}_t$ for any $t$. Let $\tau$ be geodesic in $\BH^2$ with $\PI\tau = \{0, \pi\}\subset \PI \BH^2$. Let $\alpha=P\cap (\tau\times\BR)$. Then, $\alpha$ is the generating curve for $P$ where $\PI \alpha= \{(0,c),(0,-c)\}$, i.e. $P=\bigcup_t\wh{\varphi}_t(\alpha)$. 

On the other hand, let $P_h$ be the minimal plane with $\PI P_h=\partial R_h$ where the height of the rectangle $R_h$ is $h$, i.e. $2c=h$. The invariance of $P_h$ under the isometry $\wh{\varphi}$ shows that
$\gamma_h=P_h\cap \BH^2\times\{0\}$ is an equidistant curve from the geodesic $\wh{\gamma}=\gamma\times\{0\}$ in $\BH^2\times\{0\}$. Let $d_h=dist(\gamma_h,\wh{\gamma})$. Then, they also show that if $h\nearrow \infty$ then $d_h\searrow 0$ and if $h\searrow \pi$ then $d_h\nearrow \infty$. In other words, when $h\to \infty$, $P_h$ gets closer to the vertical geodesic plane $\gamma\times \BR$. When $h\searrow \pi$, $P_h$ escapes to infinity. Moreover, the upper half of  $P_h$, $P_h\cap \BH^2\times [0,c]$, is a vertical graph over the  component of $\BH^2\times\{0\}-\gamma_h$ in the $R\subset\Si$ side.

Now, we  show that tall rectangles are indeed quite special. They bound a unique area minimizing surface which is area minimizing.

\begin{lem}  \label{rectangle} [Tall Rectangles are Strongly Fillable] If $R$ is a tall rectangle in $\Si$, then there exists a unique minimal surface $P$ in $\BHH$ with $\PI P = \partial R$. Furthermore, $P$ is also an area minimizing surface in $\BHH$.
\end{lem}

\begin{pf} Outline of the proof is as follows. By using rectangles $\wh{R}_h \subset \Si$, we  foliate a convex region $\Delta$ in $\BHH$ by minimal planes $\wh{P}_h$ with $\PI\wh{P}_h=\partial \wh{R}_h$. As our minimal plane $P=\wh{P}_{h_0}$ is a leave in this foliation, it is the unique minimal surface bounding  $\Gamma_{h_0}=\partial \wh{R}_{h_0}$, and hence area minimizing.
	
\vspace{.2cm}
	
\noindent {\bf Step 1:} Defining the convex region $\Delta$.

\vspace{.2cm}
	
The convex region $\Delta$ will be a component of the complement of a vertical geodesic plane in $\BHH$, i.e. $\BHH-(\eta\times\BR)$. The setup is as follows:
Let $R_{h}=[-\theta_1,\theta_1]\times[-h,h]$ be a tall rectangle in $\Si$, i.e. $h>\frac{\pi}{2}$ and $0<\theta_1<\pi$. By Lemma \ref{stdisk}, for any $h>\frac{\pi}{2}$, there exists a minimal surface $P_h$ with $\PI P_h = \Gamma_h=\partial R_h$. Moreover, by the construction \cite{ST}, $\{P_h\}$ is a continuous family of complete minimal planes with $P_h\cap P_{h'}=\emptyset$ for $h\neq h'$. Now, fix $h_0>\frac{\pi}{2}$, and let $R_{h_0}=[-\theta_1,\theta_1]\times[-h_0,h_0]$

Let $\tau$ be geodesic in $\BH^2$ with $\PI\tau = \{0, \pi\}\subset \PI \BH^2$. Let $\psi_t$ be the hyperbolic isometry of $\BH^2$ which fixes $\tau$, where $t$ is the translation parameter along $\tau$. In particular, in the upper half plane model $\BH^2=\{(x,y) \ | \ y>0\}$, $\tau=\{(0,y) \ | \ y>0\}$ and $\psi_t(\mathbf{x})=t\mathbf{x}$. Then, let  $\theta_t=\psi_t(\theta_1)$. Then for $0<t<\infty$, $0<\theta_t<\pi$. Hence, $\theta_t<\theta_1$ when $0<t<1$, and $\theta_t>\theta_1$ when $1<t<\infty$. In particular, this implies  $[-\theta_1,\theta_1]\subset [-\theta_t,\theta_t]$ for $t>1$, and $[-\theta_1,\theta_1]\supset [-\theta_t,\theta_t]$ for $t<1$. For notation, let $\theta_0=0$ and let $\theta_\infty=\pi$.

Now, define a continuous family of rectangles $\wh{R}_h$ which foliates an infinite vertical strip in $\Si$ as follows. Let $s:(\frac{\pi}{2}, \infty) \to (0, 2)$ be a smooth monotone increasing function such that $s(h)\nearrow 2$ when $h\nearrow \infty$, and $s(h)\searrow 0$ when $h\searrow\frac{\pi}{2} $. Furthermore, let $s(h_0)=1$.

Now, define $\wh{R}_h$ be the rectangle in $\Si$ with $\wh{R}_h=[-\theta_{s(h)},\theta_{s(h)}]\times[-h,h]$. Hence, $\wh{R}_{h_0}=R_{h_0}$, and for any $h\in (\frac{\pi}{2}, \infty)$, $\wh{R}_h$ is a tall rectangle with height $2h>\pi$. Let $\wh{\Gamma}_h=\partial \wh{R}_h$. Then, the family of simple closed curves $\{\wh{\Gamma}_h\}$ foliates the vertical infinite strip $\Omega=((-\theta_2 ,\theta_2 ) \times \BR) - (\{0\}\times [-\frac{\pi}{2},\frac{\pi}{2}])$ in $\Si$.

Recall that  $R_{h}=[-\theta_1,\theta_1]\times[-h,h]$ for any $h>\pi/2$, and the planes $P_h$ are minimal surfaces with $\PI P_h = \Gamma_h$. Let $\wh{\psi}_t$ be the isometry of $\BHH$ with $\wh{\psi}_t(p, s)=(\psi_t(p),s)$ where $p\in \BH^2$ and $s\in \BR$. Then clearly $\wh{R}_h=\wh{\psi}_{s(h)}(R_h)$. In other words, $\wh{R}_h$ and $R_h$ have the same height, but $\wh{R}_h$ is "widened $R_h$" in the horizontal direction via isometry $\wh{\psi}$. Similarly, define $\wh{P}_h=\wh{\psi}_{s(h)}(P_h)$. Hence, $\wh{P}_h$ is a complete minimal plane with $\PI \wh{P}_h=\wh{\Gamma}_h=\partial \wh{R}_h$. 

Notice that $\wh{P}_\infty$ is the geodesic plane $\eta\times\BR$ in $\BHH$ where $\eta$ is a geodesic in $\BH^2$ with $\PI \eta = \{-\theta_2,\theta_2\}$. Let $\Delta$ be the component of $\BHH-\wh{P}_\infty$ containing $P_{h_0}$, i.e. $\partial \overline{\Delta}=\wh{P}_\infty$ and $\PI \overline{\Delta} = \overline{\Omega}$. We claim that the family of complete minimal planes $\{\wh{P}_h \ | \ h\in(\frac{\pi}{2},\infty)\}$ foliates $\Delta$.

\vspace{.2cm}

\noindent {\bf Step 2:} Foliating $\Delta$ by minimal planes $\{\wh{P}_h\}$.

\vspace{.2cm}

Notice that as $\{P_h\}$ is a continuous family of minimal planes, and $\{\wh{\psi}_t\}$ is a continuous family of isometries, then by construction $\wh{P}_h=\wh{\psi}_{s(h)}(P_h)$ is a continuous family of minimal planes, and $\Delta=\bigcup_{h\in(\frac{\pi}{2},\infty)} \wh{P}_h$. Hence, all we need to show that $\wh{P}_h\cap\wh{P}_{h'}=\emptyset$ for $h<h'$. First notice that $P_h\cap P_{h'}=\emptyset$ by \cite{ST}. Hence, $\wh{\psi}_{s(h)}(P_h) \cap \wh{\psi}_{s(h)} (P_{h'}) =\emptyset$. Let $s'=s(h')/s(h)>1$.

Notice that both planes $\wh{\psi}_{s(h)}(P_h)$ and $\wh{\psi}_{s(h)} (P_{h'})$ are graphs over rectangles $[-\theta_{s(h)},\theta_{s(h)}] \times [-h,h]$ and $[-\theta_{s(h)},\theta_{s(h)}] \times [-h',h']$ respectively. For any $c\in (-h,h)$, the line $l^{h'}_c = \wh{\psi}_{s(h)} (P_{h'})\cap (\BH^2\times\{c\})$ is on far side ($\pi\in \si$ side) of the line $l^{h}_c = \wh{\psi}_{s(h)} (P_{h})\cap \BH^2\times\{c\}$ in $\BH^2\times \{c\}$. Hence, for any $c$, $\psi_{s'}(l^{h'}_c)\cap l^h_c = \emptyset$ since $\psi_{s'}$ pushes $\BH^2$ toward $\pi\in \PI\BH^2$ as $s'>1$. As $\wh{\psi}_{s'}\circ \wh{\psi}_s(h)=\wh{\psi}_{s'.s(h)}=\wh{\psi}_{s(h')}$, then $\wh{\psi}_{s(h)}(P_h) \cap \wh{\psi}_{s(h')} (P_{h'}) =\emptyset$. In other words, $\wh{P}_h\cap\wh{P}_{h'}=\emptyset$ for $h<h'$. In particular, $\{\wh{P}_h\}$ is a pairwise disjoint family of planes, with $\Delta=\bigcup_{\frac{\pi}{2}}^\infty \wh{P}_h$. This shows that the family of minimal planes $\{\wh{P}_h \ | \ h\in(\frac{\pi}{2},\infty)\}$ foliates $\Delta$.

\vspace{.2cm}

\noindent {\bf Step 3:} $P_{h_0}$ is the unique minimal surface with asymptotic boundary $\Gamma_{h_0}=\partial R_{h_0}$ in $\Si$, i.e $\PI P_{h_0}=\Gamma_{h_0}$.

\vspace{.2cm}

Assume on the contrary. If there was another minimal surface $\Sigma$ in $\BHH$ with $\PI\Sigma = \partial R_{h_0}$, then $\Sigma$ must belong to the convex region $\Delta$ by the convex hull principle. In particular, one can easily see this fact by foliating $\BHH-\Delta$ by the geodesic planes $\{\wh{\psi}_t(\wh{P}_\infty) \ | \ t>1\}$. Hence, if $\Sigma\nsubseteq \Delta$, then for $t_0=\sup_t\{\Sigma\cap \wh{\psi}_t(\wh{P}_\infty)\neq \emptyset\}$, $\Sigma$ would intersect the geodesic plane $\wh{\psi}_{t_0}(\wh{P}_\infty)$ tangentially with lying in one side. This contradicts to maximum principle as both are minimal surfaces.

Now, since $\Sigma\subset \Delta$ and $\Delta$ is foliated by $\wh{P}_h$, if $\Sigma\neq P_{h_o}$, then $\Sigma\cap P_h\neq \emptyset$ for some $h\neq h_o$. Then, either $h_1=\sup\{h>h_o\ | \ \Sigma\cap \wh{P}_h \neq \emptyset\}$ or $h'_1=\inf\{h<h_o\ | \ \Sigma\cap \wh{P}_h \neq \emptyset\}$ exists. In either case, $\Sigma$ would intersect $\wh{P}_{h_1}$ or $\wh{P}_{h_1'}$ tangentially by lying in one side. Again, this contradicts to maximum principle as both are minimal surfaces. Hence, such a $\Sigma$ cannot exist, and the uniqueness follows.

\vspace{.2cm}

\noindent {\bf Step 4:} $P_{h_0}$ is indeed an area minimizing surface in $\BHH$.

\vspace{.2cm}

Now, we  finish the proof by showing that $P_{h_0}$ is indeed an area minimizing surface in $\BHH$. Let $B_n$ be the $n$-disk in $\BH^2$ with the center origin $O$ in the Poincare disk model, i.e. $B_n=\{x\in\BH^2 \ | \ d(x,O)<n\}$. Let $\wh{B}_n=B_n\times [-h_0,h_0]$ in $\BHH$.  We claim that $P^n_{h_0}=P_{h_0}\cap \wh{B}_n$ is an area minimizing surface, i.e. $P^n_{h_0}$ has the smallest area among the surfaces $S$ in $\BHH$ with the same boundary, i.e. $\partial P^n_{h_0}=\partial S \ \Rightarrow \ |P^n_{h_0}|\leq |S|$ where $|.|$ represents the area.

Let $\Omega_n=\wh{B}_n \cap \overline{\Delta}$ be the compact, convex subset of $\BHH$. Let $\beta_n=\partial P^n_{h_0}$ be the simple closed curve in $\partial \Omega_n$. Notice that by the existence theorem of area minimizing surfaces (Lemma \ref{AMSexist}), there exists an area minimizing surface $\Sigma$ in $\BHH$ with $\partial \Sigma = \beta_n$. Furthermore, as $\Omega_n$ is convex, $\Sigma\subset \Omega_n$. However, as $\{\wh{P}_h \ | \ h\in(\frac{\pi}{2},\infty)\}$ foliates $\Delta$, $\{\wh{P}_h \cap \Omega_n\}$ foliates $\Omega_n$. Similar to above argument, if $\Sigma$ is not a leaf of this foliation, there must be a last point of contact with the leaves, which gives a contradiction with the maximum principle. Hence, $\Sigma=P^n_{h_0}$, and $P^n_{h_0}$ is an area minimizing surface. This shows that any compact subsurface of $P_{h_0}$ is an area minimizing surface as it must belong to $P^n_{h_0}$ for sufficiently large $n>0$. This proves $P_{h_0}$ is an area minimizing surface with $\PI P_{h_0}=\Gamma_{h_0}$, and it is the unique minimal surface in $\BHH$ with asymptotic boundary $\partial R_{h_0}$ in $\Si$. As any tall rectangle in $\Si$ is isometric image of $R_h$ for some $\frac{\pi}{2}<h<\infty$, the proof follows.
\end{pf}

\subsection{Tall Curves} \ 

\vspace{.2cm}

After defining, and studying tall rectangles in $\Si$ (Section \ref{tallrecsec}), now we are ready to define tall curves in $\Si$.

\begin{figure}[t]
	\begin{center}
		$\begin{array}{c@{\hspace{.4in}}c}

		\relabelbox  {\epsfysize=2in \epsfbox{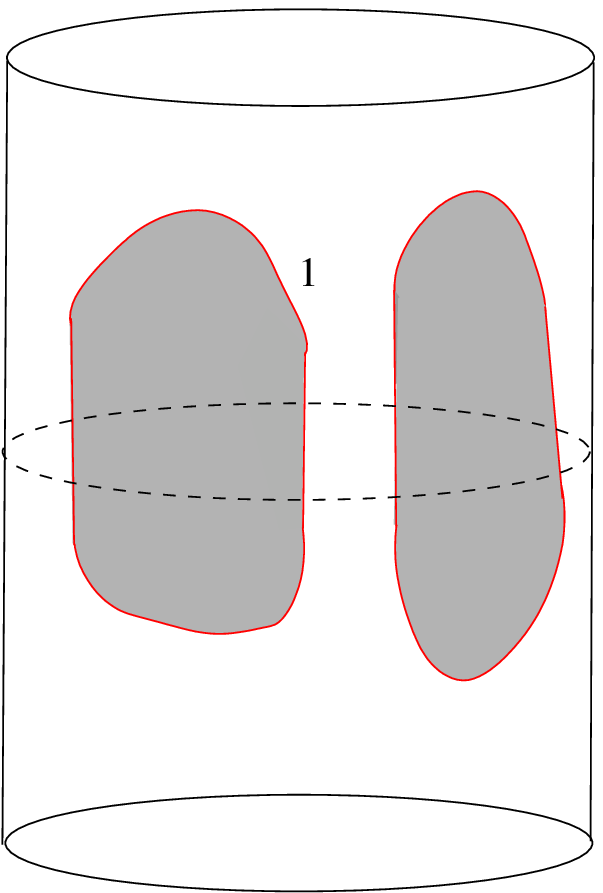}} \relabel{1}{\small $\Omega^+$} \endrelabelbox &
		
		\relabelbox  {\epsfysize=2in \epsfbox{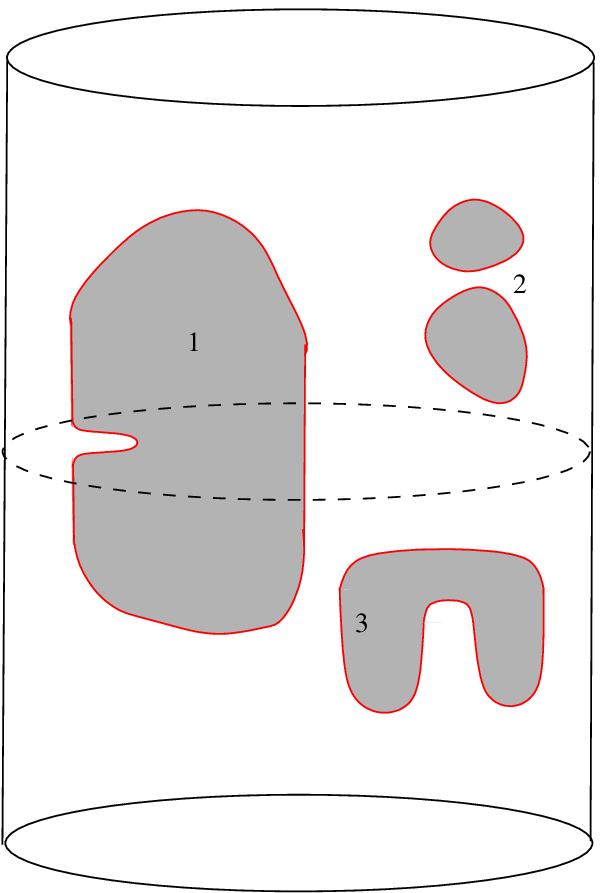}} \relabel{1}{\tiny $\Omega^-_1$} \relabel{2}{\tiny $\Omega^+_2$} \relabel{3}{\tiny $\Omega^-_3$}   \endrelabelbox \\
	\end{array}$
	
\end{center}
\caption{\label{tall} \footnotesize In the left, $\Gamma$ is a tall curve with two components. In the right, there are three nonexamples of tall curves. Shaded regions describe the $\Omega^-_i$ where $\Si-\Gamma_i=\Omega^+_i\cup\Omega^-_i$.}
\end{figure}

\begin{defn} \label{talldef} [Tall Curves] 	We call a finite collection of disjoint simple closed curves $\Gamma$ in $\Si$  {\em tall curve} if the region $\Gamma^c=\Si-\Gamma$ can be written as a union of open tall rectangles $R_i=(\theta^i_1,\theta^i_2)\times(t^i_1,t^i_2)$, i.e. $\Gamma^c=\bigcup_i R_i$ (See Figure \ref{tall}).
	
We  call a region $\Omega$ in $\Si$ a {\em tall region}, if $\Omega$ can be written as a union of tall rectangles, i.e. $\Omega=\bigcup_i R_i$ where  $R_i$ is a tall rectangle.
\end{defn}

On the other hand, by using the idea above, we can define a notion called {\em height of a curve} as follows:

\begin{defn} \label{height} [Height of a Curve] Let $\Gamma$ be a collection of simple closed curves in $\Si$, and let $\Omega=\Si-\Gamma$. For any $\theta\in[0,2\pi)$, let $L_\theta=\{\theta\}\times \BR$ be the vertical line in $\Si$. Let $L_\theta\cap\Omega=l_\theta^1\cup..\cup l_\theta^{i_\theta}$ where $l_\theta^i$ is a component of $L_\theta\cap \Omega$. Define the height $h(\Gamma)= \inf_\theta\{|l_\theta^i|\}$.

Notice that $\Gamma$ is a tall curve if and only if $h(\Gamma)>\pi$. Now, we say $\Gamma$ is a {\em short curve} if $h(\Gamma)<\pi$.
\end{defn}

\begin{rmk} Note that if $\Gamma$ is a finite collection of disjoint simple closed curves in $\Si$, then we can always write $\Gamma^c=\Omega^+\cup\Omega^-$ where $\Omega^\pm$ are (possibly disconnected) tall regions with $\partial\overline{\Omega^+}=\partial\overline{\Omega^-}=\Gamma$. Notice that if $\Gamma$ has more than one component, than $\Omega^+$ or $\Omega^-$ may not be connected.

Note also that any curve containing a thin tail is short curve by definition. However, there are some short curves with no thin tails, like $\Gamma_3$ in Figure \ref{tall}-right and Figure \ref{barrierfig}-right.

Notice also that for each nullhomotopic component $\gamma_i$ of a tall curve $\Gamma$, if  $\theta_i^+$ ($\theta_i^-$) is a local maximum (minimum) of horizontal coordinates of $\gamma_i$, then by Lemma \ref{nonexist}, $L_{\theta_i^\pm}\cap\gamma_i$ must be a pair of vertical line segments of length greater than $\pi$ (See Figure \ref{tall} left). Also, in Figure \ref{tall} right, three non-tall curves $\Gamma_1, \Gamma_2$ and $\Gamma_3$ are pictured as examples. If we name the shaded regions as $\Omega_i^-$, $\Gamma_1$ is not tall as $\Omega^+_1$ is not tall because of the small cove. $\Gamma_2$ has two components, and it is not tall as $\Omega_2^+$ is not tall (The two components are very close to each other). Finally, $\Gamma_3$ is not tall as $\Omega^-_3$ is not tall region because of the short neck.
\end{rmk}

Note also that recently in \cite{KMR}, Klaser, Menezes and Ramos generalized the {\em tall curve}, and {\em height of a curve} notions to the other $\mathbb{E}(-1,\tau)$ homogeneous spaces, and obtained several existence and nonexistence results for the asymptotic Plateau problem in these spaces.

\begin{rmk} \label{exceptional} (Exceptional Curves) We  call a \underline{short} curve $\Gamma$ {\em exceptional} if $\overline{\Omega}^\pm$ can be written as a union of closed tall rectangles ($[\theta_1,\theta_2]\times[t_1,t_2]$). Notice that as $\Gamma$ is a short curve, $\Omega^\pm$ cannot be written as a union of open tall rectangles ($(\theta_1,\theta_2)\times(t_1,t_2)$). 
\end{rmk}

\begin{wrapfigure}{r}{1in}
\vspace{-.3cm}
\relabelbox  {\epsfxsize=1in
	
	\centerline{\epsfbox{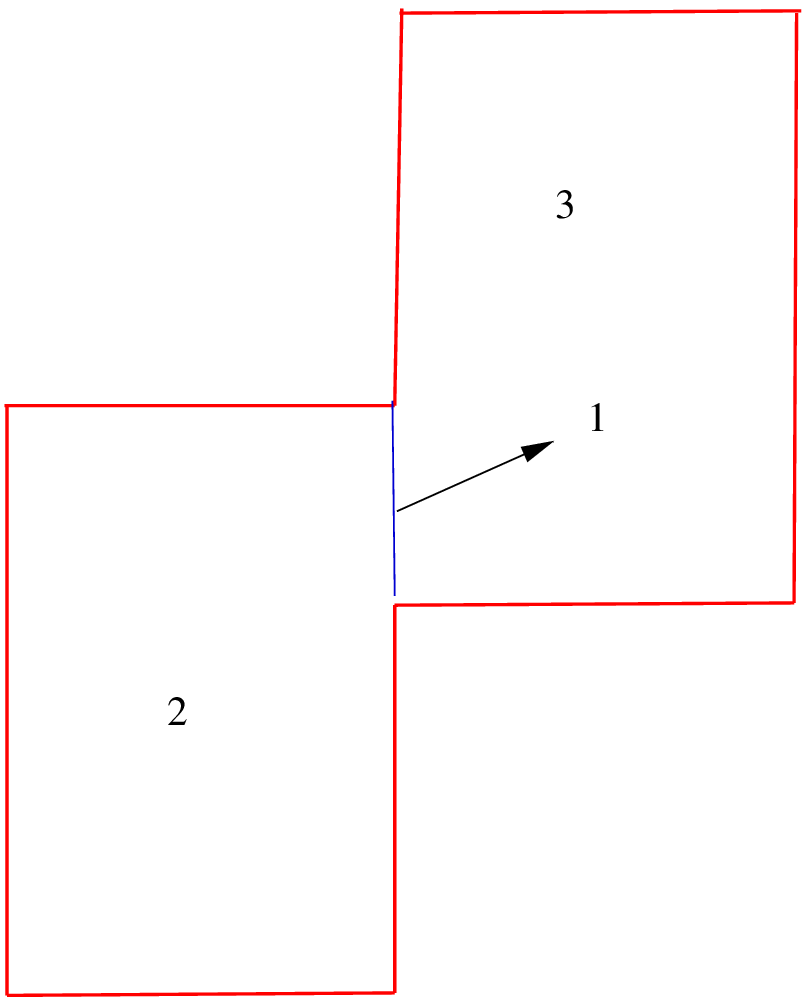}}}
\relabel{1}{\footnotesize $\mu$}
\relabel{2}{\footnotesize $R_1$}
\relabel{3}{\footnotesize $R_2$}
\endrelabelbox
\vspace{-.3cm}
\end{wrapfigure}

As an example, consider $R_1=[0,\frac{\pi}{3}]\times[-1,5]$ and $R_2=[\frac{\pi}{3}, \frac{2\pi}{3}]\times[-5,1]$. Let $\gamma_i=\partial R_i$. Let $\mu=\{\frac{\pi}{3}\}\times (-1,1)$ be a line segment of length $2$ (See Figure right). Define $\Gamma=\gamma_1\cup\gamma_2 - l$. Clearly, $h(\Gamma)=d((\frac{\pi}{3},1),(\frac{\pi}{3},-1))=\|\mu\|=2$ and $\Gamma$ is a short curve. However, $\overline{\Omega}^+=\overline{R}_1\cup\overline{R}_2$ which makes $\Gamma$ an exceptional curve. These curves are very small set of curves among the space of curves in $\Si$, however they will have a problematic feature with barrier argument when we show the {\em nonexistence} for short curves in Theorem \ref{APP}. So, throughout the paper, we will assume that closed curves in $\Si$ are not exceptional unless otherwise stated. We would like to thank Laurent Mazet for pointing out the exceptional curves.

\section{Asymptotic Plateau Problem in $\BHH$} \label{secAPP}

In this section, we prove our main result. Note that for the following theorem, we assume $\Gamma$ is not an exceptional curve (See Remark \ref{exceptional}).

\begin{thm} \label{APP} [Asymptotic Plateau Problem for $\BHH$] Let $\Gamma$ be a finite collection of disjoint Jordan curves in $\Si$ with $h(\Gamma)\neq \pi$. Then, there exists an  area minimizing surface $\Sigma$ in $\BHH$ with $\PI \Sigma = \Gamma$ if and only if $\Gamma$ is a tall curve. Furthermore, all such surfaces are embedded.
\end{thm}

\noindent {\em Outline of the proof:} We  use the standard techniques for the asymptotic Plateau problem \cite{An}. In particular, we  construct a sequence of compact area minimizing surfaces $\{\Sigma_n\}$ in $\BHH$ with $\partial \Sigma_n\to \Gamma$, and in the limit, we aim to obtain an area minimizing surface $\Sigma$ with $\PI \Sigma=\Gamma$. Notice that the main issue here is not to show that $\Sigma$ is an area minimizing surface, but to show that $\Sigma$ is not escaping to infinity, i.e. $\Sigma \neq \emptyset$ and $\PI \Sigma = \Gamma$ (See Remark \ref{escape2}). Recall that by Lemma \ref{nonexist}, if a simple closed curve $\gamma$ in $\Si$ has a thin tail, then there is no minimal surface $S$ in $\BHH$ with $\PI S = \gamma$. This means that if you similarly construct area minimizing surfaces $S_n$ with $\partial S_n\to \gamma$, then either $S=\lim S_n = \emptyset$ or $\PI S\neq \gamma$, i.e. the sequence $S_n$ escapes to infinity completely ($S=\emptyset$), or some parts of the sequence $S_n$ escapes to infinity ($\PI S \neq \gamma$). 

In particular, in the following, we aim to show that for a tall curve $\Gamma$, the limit surface $\Sigma$ does not escape to infinity, and $\PI\Sigma=\Gamma$. We  achieve this, by constructing barriers near infinity preventing escaping to infinity.

\vspace{.2cm}

\begin{pf} We split the proof into two parts. In the first part, we show the "if" part. In the second part, we  prove the converse.
	
\vspace{.2cm}

\noindent {\bf Step 1:} [Existence] If $\Gamma$ is tall ($h(\Gamma)>\pi$), then there exists an area minimizing surface $\Sigma$ in $\BHH$ with $\PI \Sigma = \Gamma$.

\vspace{.2cm}

\noindent {\bf Step 1A:} Construction of the barrier $\X$ near $\Si$.

\vspace{.2cm}

\noindent {\em Proof of Step 1A:} In this part, we  construct a barrier $\X$ which prevents the limit escape to infinity.
	
Since $\Gamma$ is a tall curve, by definition, $\Gamma^c=\Si-\Gamma= \Omega^+\cup\Omega^-$ where $\Omega^\pm$ is a tall open region with $\partial \overline{\Omega^\pm}=\Gamma$. Notice that if $\Gamma$ has more than one component, $\Omega^+$ or $\Omega^-$ may not be connected. Note also that for any component of $\Gamma$, one side belongs to $\Omega^+$, and the other side belongs to $\Omega^-$ by assumption.

Let $\Omega^\pm=\bigcup_{\alpha\in\A^\pm} R_\alpha$ where $\{R_\alpha\}$ are tall rectangles in $\Si$. For each tall rectangle $R_\alpha$, by Lemma \ref{rectangle}, there exists a unique area minimizing surface $P_\alpha$ with $\PI P_\alpha = \partial R_\alpha$. Let $\Delta_\alpha$ is the component of $\BHH-P_\alpha$ with $\PI\Delta_\alpha=R_\alpha$. Define $\X^\pm=\bigcup_{\alpha\in\A^\pm}\Delta_\alpha$. Then, by construction $\PI\X^\pm=\Omega^\pm$. Let $\X=\X^+\cup\X^-$. We call $\X$ is a barrier near infinity. Notice that $\X$ is an open region in $\BHH$ with $\PI\X=\Si-\Gamma$.

\vspace{.2cm}

\noindent {\bf Step 1B:} Construction of the sequence $\{\Sigma_n\}$.

\vspace{.2cm}

\noindent {\em Proof of Step 1B:} Let $C>0$ be sufficiently large that $\Gamma\subset \PI \BH^2 \times (-C,C)$. Let $B_n$ be the $n$-disk in $\BH^2$ with the center origin, and $\wh{B}_n = B_n\times [-C,C]$ is an compact solid cylinder in $\BHH$. Let $\Gamma_n$ be the radial projection of $\Gamma$ into the cylinder $\partial B_n\times [-C,C]$. Then, $\Gamma_n$ is a finite union of disjoint Jordan curves in $\partial \wh{B}_n$. Notice that for any $\alpha\in \A^\pm$, $P_\alpha$ is a graph over $R_\alpha$ by Section \ref{tallcurvesec}. This implies for any $n$, $\Gamma_n\cap\X=\emptyset$ by the construction of $\X$. 

Let $\Sigma_n$ be the area minimizing surface in $\BHH$ with $\partial \Sigma_n = \Gamma_n$ by Lemma \ref{AMSexist}. Then, as $\wh{B}_n$ is convex, $\Sigma_n\subset \wh{B}_n$.

\vspace{.2cm}

\noindent {\bf Step 1C:} For any $n$, $\Sigma_n\cap \X=\emptyset$.

\vspace{.2cm}

\noindent {\em Proof of Step 1C:}  Recall that $\X^\pm=\bigcup_{\alpha\in\A^\pm} \Delta_\alpha$. Hence, we can show that for any $\alpha\in\A^\pm$, $\Sigma_n\cap P_\alpha=\emptyset$, we are done. Fix $\alpha_0\in \A^+$. 
Let $R_{\alpha_o}$ be a tall rectangle with $R_{\alpha_o}\subset \Omega^+$. Let $P_{\alpha_o}$ be the unique area minimizing surface with $\PI P_{\alpha_o}=\partial R_{\alpha_o}$. We claim that $\Sigma_n\cap P_{\alpha_o}=\emptyset$ for any $n$. Let $P^n_{\alpha_o} = P_{\alpha_o}\cap \wh{B}_n$. Recall that $P_\alpha$ is a graph over $R_\alpha$, and $\Gamma_n$ is radial projection of $\Gamma$ to $\wh{B}_n$. Therefore, $\partial P^n_{\alpha_o}= \eta_n$ and $\partial \Sigma_n=\Gamma_n$ are disjoint simple closed curves in $\partial \wh{B}_n$. Assume that $\Sigma_n\cap P^n_{\alpha_o}\neq \emptyset$. Then, as both are separating in $\wh{B}_n$, the intersection must consist of a collection of closed curves $\{\mu_1,..\mu_k\}$ (no isolated points in the intersection because of the maximum principle). Let $T_n$ be a component of $\Sigma_n-P^n_{\alpha_o}$ with $\Gamma_n\nsubseteq \partial T_n$. Let $Q_n\subset P^n_{\alpha_o}$ be the collection of disks with $\partial Q_n = \partial T_n$. Since both $\Sigma_n$ and $P^n_{\alpha_o}$ are area minimizing, then so are $T_n$ and $Q_n$. Hence, they have the same area $|T_n|=|Q_n|$ as $\partial T_n = \partial Q_n$. Let $\Sigma_n'=(\Sigma_n-T_n)\cup Q_n$. Then, since $\partial \Sigma_n = \partial \Sigma_n'$ and $|\Sigma_n|=|\Sigma_n'|$, $\Sigma_n'$ is also area minimizing surface. However, $ \Sigma_n'$
is not smooth along $\partial Q_n$ which contradicts to the interior regularity of area minimizing surfaces (Lemma \ref{AMSexist}). This shows that $\Sigma_n\cap P_{\alpha_o}=\emptyset$ for any $n$. Hence,  Step 1C follows.

\vspace{.2cm}

\noindent {\bf Step 1D:} The limit area minimizing surface $\Sigma$ is not empty.

\vspace{.2cm}

\noindent {\em Proof of Step 1D:} As described in the outline at the beginning of the proof, first we need to guarantee that the sequence $\{\Sigma_n\}$ is not escaping to infinity, i.e. $\lim\Sigma_n=\Sigma\neq \emptyset$. Let $\Sigma$ be the limit of $\Sigma_n$. In particular, by the convergence theorem (Lemma \ref{convergence}), for any compact solid cylinder $\wh{B}_m$, the sequence $\{\Sigma_n\cap \wh{B}_m\}$ has a convergent subsequence with limit $\Sigma^m\subset \wh{B}_m$. By using the diagonal sequence argument, in the limit, we get an area minimizing surface $\Sigma$ with $\Sigma\cap \wh{B}_m=\Sigma^m$. Notice also that $\Sigma^m$ separates $\wh{B}_m$ where the component near boundary contains $P^m_{\alpha_o}$ as $\Sigma_n\cap P_{\alpha_o}=\emptyset$ for any $n$. Hence, if $P_{\alpha_o}\cap \wh{B}_m \neq \emptyset$, then $\Sigma\cap \wh{B}_m=\Sigma^m\neq \emptyset$ as $\Sigma^m$ separates $P^m_{\alpha_o}$ in $\wh{B}_m$. This also implies $\Sigma$ is not empty. In particular, for any $n$, $\Sigma_n$ stays in one side (far side from infinity) of $P_{\alpha_o}$, and $P_{\alpha_o}$ acts as a barrier which prevents the sequence $\{\Sigma_n\}$ escaping to infinity.

\vspace{.2cm}

\noindent {\bf Step 1E:} $\PI\Sigma=\Gamma$.

\vspace{.2cm}

\noindent {\em Proof of Step 1E:} First, we show that $\PI\Sigma \subset \Gamma$. By Step 1C, $\Sigma_n\cap \X=\emptyset$. As $\X$ is open, this implies $\Sigma\cap \X=\emptyset$. As $\PI X=\Si-\Gamma$, we have $\PI\Sigma\subset \Gamma$.


We finish the proof by showing that $\PI\Sigma\supset\Gamma$. Let $p\in \Gamma$. We will show that $p\in \overline{\Sigma}$. Let $p$ be in the component $\gamma$ in $\Gamma$. As $\Gamma^c=\Omega^+\cup\Omega^-$, let $\{p_i^\pm\}\subset \Omega^\pm$ be two sequences in opposite sides of $\gamma$ with $\lim p_i^\pm=p$. Let $\alpha_i$ be a small circular arc in $\overline{\BHH}$ with $\partial \alpha_i = \{p_i^+,p_i^-\}$ and $\alpha_i \bot \Si$. Then, for any $i$, there exists $N_i$ such that for any $n>N_i$, $\Gamma_n$ links $\alpha_i$, i.e. $\Gamma_n$ is not nullhomologous in $\BHH-\alpha_i$. Hence, for any $n>N_i$, $\alpha_i\cap \Sigma_n\neq \emptyset$. This implies $\Sigma \cap\alpha_i\neq \emptyset$ for any $i$ by construction. Like above, let $R^\pm_{i}\subset \Omega^\pm$ be the tall rectangle with $p^\pm_i\subset R^\pm_{i}$. Similarly, let $P^\pm_{i}$ be the unique area minimizing surface with $\PI P^\pm_{i}=\partial R^\pm_{i}$. Let $\alpha'_i\subset \alpha_i$ be a subarc with $\partial \alpha'_i \subset P^+_{i}\cup P^-_{i}$. Hence, $\alpha'_i$ is a compact arc in $\BHH$. Moreover, as $\Sigma\cap P^\pm_{i}=\emptyset$, then there exists a point $x_i$ in $\Sigma \cap\alpha_i'$ for any $i$. Then, $\lim  x_i=p$. Hence $p\in \overline{\Sigma}$, and  $\PI\Sigma=\Gamma$. Step 1 follows.

\vspace{.2cm}

\noindent {\bf Step 2:} [Nonexistence] If $\Gamma$ is short ($h(\Gamma)<\pi$), then there is no area minimizing surface $\Sigma$ in $\BHH$ with $\PI \Sigma = \Gamma$.

\begin{pf} Assume that there exists an area minimizing surface $\Sigma$ in $\BHH$ with $\PI\Sigma = \Gamma$. Note that we a priori assume that $\Gamma$ is not an exceptional curve (See Remark \ref{exceptional}). Since $\Gamma$ is a short curve in $\Si$, there is a $\theta_0\in[0,2\pi)$ with $(\theta_0,c_1),(\theta_0,c_2)\in \Gamma$ where $0<c_1-c_2<\pi-2\e$ for some $\e>0$. Let $c^+=c_1+\e$ and $c^-=c_2-\e$. Let $p^+=(\theta_0,c^+)$ and $p^-=(\theta_0,c^-)$ where $p^\pm\not\in \Gamma$. 

Since $\PI\Sigma=\Gamma$, this implies $\overline{\Sigma}= \Sigma\cup\Gamma$ is a surface with boundary in $\overline{\BHH}$ by Lemma \ref{surface}. Let $d_E$ be the Euclidean metric on $\overline{\BHH}$, define $O^\pm=\{q\in\overline{\BHH} \ | \ d_E(q,p^\pm)<\delta_1\}$ as an open neighborhood of $p^\pm$ in $\overline{\BHH}$ such that $O^\pm\cap\overline{\Sigma}=\emptyset$.  Let $D^\pm=(\BH^2\times \{c^\pm\} )\cap O^\pm$. By construction, $D^\pm$ contains a half plane in the hyperbolic plane $\BH^2\times \{c^\pm\}$.

By Lemma \ref{cat1} and Lemma \ref{cat2} in the Appendix, for any $h<\pi$, there exist an {\em area minimizing} compact catenoid $S$ of height $h$. For $h=c^+-c^-<\pi$, let $S$ be the area minimizing compact catenoid with  $\partial S\subset \BH^2\times\{c^-,c^+\}$. In other words, $\partial S$ consists of two curves $\gamma^+$ and $\gamma^-$ where $\gamma^\pm$ is a round circle of radius $\wh{\rho}(d)$ in $\BH^2\times\{c^\pm\}$ centered at the origin. Let $\theta_1$ be the antipodal point of $\theta_0$ in $\si$. Let $\psi_t$ be the hyperbolic isometry fixing the geodesic between $\theta_0$ and $\theta_1$.  In particular, $\psi_t$ corresponds to $\psi_t(x,y)=(tx,ty)$ in the upper half space model where $\theta_1$ corresponds to origin, and $\theta_0$ corresponds to the point at infinity. Let $\wh{\psi}_t:\BHH\to\BHH$ be the isometry of $\BHH$ where $\wh{\psi}_t(p,z)=(\psi_t(p),z)$.

Let $S_t=\wh{\psi}_t(S)$ be the isometric image of the area minimizing catenoid $S$ in $\BHH$. Let $\partial S_t =\gamma^+_t\cup\gamma^-_t$ where $\gamma^\pm_t=\psi_t(\gamma^\pm)$. Notice that $\gamma^\pm_t\subset\BH^2\times\{c^\pm\}$. 
Let $N_o>0$ be sufficiently large that $\gamma^+_t\subset D^+$ and $\gamma^-_t\subset D^-$ for any $t>N_o$. Then, for any $t>N_o$, $\partial S_t\subset D^+\cup D^-$, and $\partial S_t \cap \Sigma=\emptyset$.

Recall $\SI-\Gamma=\Omega^+\cup \Omega^-$ where $\partial \Omega_1=\partial \Omega_2=\Gamma$. Let $\overline{\BHH}-\overline{\Sigma}= \Delta^+\cup\Delta^-$ where $\PI \Delta^\pm=\Omega^\pm$. Let $\beta=\{\theta_0\}\times(c_1,c_2)$ be the vertical line segment in $\Si$, and let $\beta\subset \Omega^+$. Since $\Delta^+$ is an open subset in $\overline{\BHH}$ and $\beta\subset \Delta_1$, then an open neighborhood $O_\beta$ of $\beta$ in $\BHH$ must belong to $\Delta^+$. Then, by construction, we can choose $t_o>N_o$ sufficiently large that $S_{t_o}\cap O_\beta\neq \emptyset$ and $S_{t_o}\cap O_\beta$ is connected. This shows that $S_{t_o}\cap\Sigma\neq \emptyset$. Let $S_{t_o}\cap\Sigma=\alpha$. Notice that as both $\Sigma$ and $S_{t_o}$ are area minimizing surfaces and $\partial S_t\cap \Sigma=\emptyset$, $\alpha$ is a collection of closed curves, and contains no isolated points because of the maximum principle.

Let $E$ be the compact subsurface of $\Sigma$ with $\partial E=\alpha$. In other words, $S_{t_o}$ separates $E$ from $\Sigma$. Similarly, let $T$ be the subsurface of $S_{t_o}$ with $\partial T=\alpha$. In particular, $T=S_{t_o}\cap \overline{\Delta^+}$. Since $S_{t_o}$ and $\Sigma$ are area minimizing surfaces, so are $T$ and $E$. As $\partial T=\partial E=\alpha$, and both are area minimizing surfaces, both have the same area, i.e. $|E|=|T|$.

Let $S'=(S_{t_o}-T)\cup E$. Then, clearly $\partial S_{t_o}=\partial S'$ and $|S_{t_o}|=|S'|$. Hence, as $S_{t_o}$ is an area minimizing surface, so is $S'$. However, $S'$ has singularity along $\alpha$. This contradicts to the regularity of area minimizing surfaces (Lemma \ref{AMSexist}). Step 2 follows.
\end{pf} \end{pf}

\begin{rmk} \label{pi-curves} [$h(\Gamma)=\pi$ case] Notice that the theorem finishes off the asymptotic Plateau problem for $\BHH$ except the case $h(\Gamma)=\pi$. Note that this case is delicate as there are strongly fillable and strongly non-fillable curves of height $\pi$. For example, if $\Gamma_1$ is a rectangle in $\Si$ with height $\pi$, then the discussion in Remark \ref{escape2} shows that $\Gamma_1$ bounds no minimal surface, hence such a $\Gamma$ is nonfillable. On the other hand, in Theorem \ref{minexist}, if we take $h_0=\pi$ and use the parabolic catenoid (\cite{Da}), it is not hard to show that the constructed surface is also area minimizing in $\BHH$ since the parabolic catenoid is also area minimizing (See Figure \ref{barrierfig}-right). These two examples show that the case $h(\Gamma)=\pi$ is very delicate. Note also that Sa Earp and Toubiana studied a relevant problem in \cite[Cor. 2.1]{ST}.
\end{rmk}

\begin{rmk} [Minimal vs. Area Minimizing] Notice that the theorem above does not say that {\em If $\gamma$ is a short curve, then there is no \underline{minimal} surface $S$ in $\BHH$ with $\PI S = \gamma$}. There are many examples of complete embedded minimal surfaces $S$ in $\BHH$ where the asymptotic boundary $\gamma$ is a short curve (e.g. butterfly curves). We postpone this question to Section \ref{secAPPmin} to discuss in detail.
	
\end{rmk}

\subsection{Convex Hull Property for Tall Curves.} \

\vspace{.2cm}

In this part, we give a natural generalization of convex hull property for asymptotic Plateau problem in $\BHH$.

\begin{defn} [Mean Convex Hull] Let $\Gamma$ be a tall curve in $\Si$. Consider the barrier $\X$ constructed in Step 1A in the proof of Theorem \ref{APP}. Define the mean convex hull of $\Gamma$ as $MCH(\Gamma)= \BHH-\X$. Notice that $MCH(\Gamma)$ is mean convex region in $\BHH$ by construction. Furthermore, $\PI MCH(\Gamma)=\Gamma$. 
\end{defn}

Analogous to convex hull property in $\BH^3$, we have the following property in $\BHH$.

\begin{cor} \label{MCHcor} [Convex Hull Property] Let $\Gamma$ be a tall curve in $\Si$. Let $S$ be a complete, embedded minimal surface in $\BHH$ with $\PI S=\Gamma$. Then, $S\subset MCH(\Gamma)$.	
\end{cor}

\begin{pf} 
The proof is similar to the convex hull property in other homogeneous ambient spaces. 
We use the same notation of the proof of Theorem \ref{APP}. In that proof, we proved that for our special sequence $\{\Sigma_n\}$ and the limit $\Sigma$, $\Sigma\cap \X=\emptyset$. However, the same proof works for any area minimizing surface $S$ with $\PI S=\Gamma$.

Recall that $\Gamma^c=\bigcup R_\alpha$ where $R_\alpha$ are tall rectangles. Let $P_\alpha$ be the unique area minimizing surfaces in $\BHH$ with $\PI P_\alpha=\partial R_\alpha$. Let $\Delta_\alpha$ be the components of $\BHH-P_\alpha$ with $\PI\Delta_\alpha = int(R_\alpha)$. Then, $\X=\bigcup_\alpha \Delta_\alpha$


Assume $S \nsubseteq MCH(\Gamma)=\X^c$. Then, $S \cap \Delta_\alpha \neq \emptyset$ for some $\alpha$. However, by the proof of Lemma \ref{rectangle}, we know that $\Delta_\alpha$ is foliated by minimal surfaces $\{P_t\mid t\in[0,\infty) \}$. Let $t_0=\sup_t\{P_t\cap\Sigma\neq \emptyset\}$. Again, by maximum principle, this is a contradiction as both $\Sigma$ and $P_{t_0}$ are minimal surfaces. The proof follows.
\end{pf}

One can visualize visualize $MCH(\Gamma)$ as follows: Assume $\Gamma\subset \si\times [c_1,c_2]$ for smallest $[c_1,c_2]$ possible. Then, $MCH(\Gamma)$ is the region in $\BH^2\times [c_1,c_2]$ where we carve out all $\Delta_\alpha$ defined by rectangles $R_\alpha \subset \Gamma^c$.

\section{Asymptotic Plateau Problem for Minimal Surfaces} \label{secAPPmin}

So far, we only dealt with the strong fillability question, i.e. detecting curves in $\Si$ bounding \textit{area minimizing surfaces} in $\BHH$. If we relax the question from "strong fillability" to only "fillability", the picture completely changes. In other words, we will see that detecting curves in $\Si$ bounding embedded \textit{minimal surfaces} is much more complicated than detecting the curves bounding embedded \textit{area minimizing surfaces}. In Theorem \ref{APP}, we gave a fairly complete answer to asymptotic Plateau problem in the strong fillability case. In this section, we will see that the classification of fillable curves is highly different.

A simple example to show the drastic change in the problem is the following: Let $\Gamma=\gamma_{1}\cup\gamma_{2}$ where $\gamma_i=\si\times\{c_i\}$ and $|c_1-c_2|<\pi$. Then clearly, $\Gamma$ is  a short curve and it bounds a complete minimal catenoid $\C_d$ by \cite{NSST} (See also appendix for further discussion on catenoids). On the other hand, the pair of geodesic planes, $\BH^2\times\{c_1\}\cup \BH^2\times\{c_2\}$, also bounds $\Gamma=\gamma_{1}\cup\gamma_{2}$. However, there is no area minimizing surface $\Sigma$ with $\PI \Sigma = \gamma_1\cup \gamma_2$ by Theorem \ref{APP}. This means neither catenoid, nor pair of geodesic planes are area minimizing, but just minimal surfaces. Hence, the following version of the problem becomes very interesting.

\vspace{.2cm}

\noindent {\bf Asymptotic Plateau Problem for Minimal Surfaces in $\BHH$:}

\noindent {\em For which curves $\Gamma$ in $\Si$, there exists an embedded minimal surface $S$ in $\BHH$ with $\PI S=\Gamma$.}

\vspace{.2cm}

In other words, \textit{which curves in $\Si$ are fillable?} Note that here we only discuss the finite curve case ($\Gamma\subset \Si$). For infinite curves case for the same question, see \cite{Co2}.

Recall that by Lemma \ref{nonexist}, for any short curve $\gamma$ in $\Si$ containing  a thin tail, there is no complete minimal surface $S$ in $\BHH$ with $\PI S = \gamma$. So, this result suggest that the minimal surface case is similar to the area minimizing surface case.

On the other hand, unlike the area minimizing surface case, it is quite easy to construct short curves with more than one component, bounding minimal surfaces in $\BHH$. Let $\Gamma=\gamma_1\cup..\cup\gamma_n$ be a finite collection of disjoint tall curves $\gamma_i$. Even though every component $\gamma_i$ is tall, because of the vertical distances between the components $\gamma_i$ and $\gamma_j$, the height $h(\Gamma)$ can be very small. So, $\Gamma$ itself might be a short curve, even though every component is a tall curve. For each component $\gamma_i$, our existence theorem (Theorem \ref{APP}) already gives an area minimizing surface $\Sigma_i$ with $\PI \Sigma_i=\gamma_i$. Hence, the surface $\wh{S}=\Sigma_1\cup..\Sigma_n$ is automatically a minimal surface with $\PI \wh{S}=\Gamma$. By using this idea, for any height $h_0>0$, we can trivially produce short curves $\Gamma$ with height $h(\Gamma)=h_0$ by choosing the components sufficiently close. e.g. the pair of horizontal geodesic planes $\BH^2\times\{c_1\}\cup \BH^2\times\{c_2\}$ with $|c_1-c_2|=h_0$.

Naturally, next question would be what if $\Gamma$ has only one component. Does $\Gamma$ need to be a tall curve to bound a minimal surface in $\BHH$? The answer is again no. Now, we  also construct {\em simple closed short curves} which bounds complete minimal surfaces in $\BHH$. The following result with the observation above shows that the asymptotic Plateau problem for minimal surfaces is very different from the  asymptotic Plateau problem for area minimizing surfaces. 


\begin{thm} \label{minexist} For any $h_0>0$, there exists a nullhomotopic simple closed curve $\Gamma$ with height $h(\Gamma)=h_0$ such that there exists a minimal surface $S$ in $\BHH$ with $\PI S = \Gamma$.
\end{thm}

\begin{pf} For $h_0>\pi$, we have tall rectangles with height $h_0$. So, we  assume $0<h_0\leq \pi$. Consider the rectangles $R^+=[s,\frac{\pi}{2}]\times [-m,m]$ and $R^-=[-\frac{\pi}{2},-s]\times [-m,m]$ where $s>0$ sufficiently small, and $m>>0$ sufficiently large, which will be fixed later. Consider another rectangle $Q=[-s,s]\times [0,h_0]$. Consider the area minimizing surfaces $P^+$ and $P^-$ with $\PI P^\pm = \partial R^\pm$. Let $\Gamma= (\partial R^+ \cup \partial R^-) \triangle \partial Q$ where $\triangle$ represents symmetric difference (See Figure \ref{barrierfig}). Notice that $h(\Gamma)=h_0$. We claim that there exists a complete embedded minimal surface $S$ in $\BHH$ with $\PI S=\Gamma$.

\begin{figure}[t]
\begin{center}
$\begin{array}{c@{\hspace{.4in}}c}

\relabelbox  {\epsfysize=1.5in \epsfbox{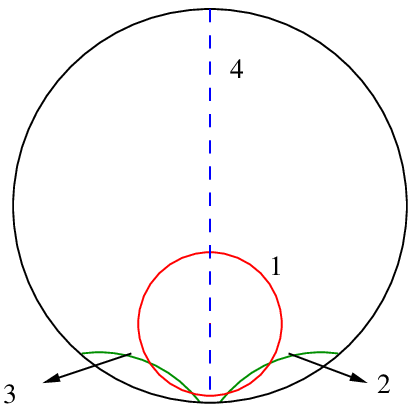}} \relabel{1}{\tiny \color{red} $C^t_{h_0}$} \relabel{2}{\tiny $P^+$} \relabel{3}{\tiny  $P^-$} \relabel{4}{\tiny \color{blue} $l$}  \endrelabelbox &

\relabelbox  {\epsfysize=2in \epsfbox{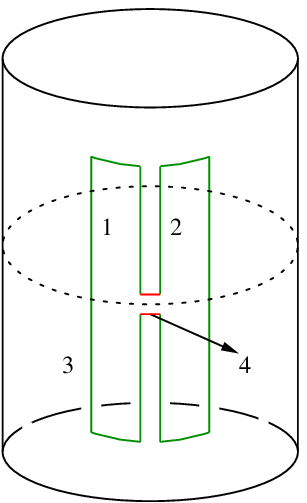}} \relabel{1}{\tiny $R^+$} \relabel{2}{\tiny $R^-$} \relabel{3}{\tiny $\Gamma$} \relabel{4}{\tiny $C^t_{h_0}$} \endrelabelbox \\
\end{array}$

\end{center}
\caption{\label{barrierfig} \footnotesize In the left, the horizontal slice $\BH^2\times \{\frac{h_0}{2}\}$ is given. In the right, $\Gamma\subset \Si$ is pictured.}
\end{figure}

Consider the minimal catenoid $C_{h_0}$ with asymptotic boundary $\si\times\{0\}\cup \si\times\{h_0\}$ (If $h_0=\pi$, take $C_{h_0}$ to be the Daniel's parabolic catenoid). Let $\varphi_t$ be the isometry of $\BHH$ which keeps $\BR$ coordinates same, fixes the geodesic $l$ in $\BH^2$ with $\PI l = \{0,\pi\}$ with translation length $\log{t}$. In particular, $\varphi_t |_{\BH^2\times \{c\}}: \BH^2\times \{c\}\to \BH^2\times \{c\}$. Furthermore, for any $p\in \BH^2\times \{c\}$, $\varphi_t(p)\to (0,c)\in \Si$ as $t\searrow 0$ and $\varphi_t(p)\to (\pi,c)\in \Si$ as $t\searrow \infty$. Now, we can choose $t>0$ and $s>0$ sufficiently small, $m>0$ sufficiently large  so that $P^+\cup P^-$ separates $\varphi_t(C_{h_0})=C^t_{h_0}$ into $4$ disks (See Figure \ref{barrierfig}). In other words, there is a component $\Delta$ in $\BHH - (P^+\cup P^-\cup C^t_{h_0})$ such that  $\PI \Delta = Q$.

Now, let $\Omega^+$ be the component of $\BHH-P^+$ such that $\PI \Omega^+ = R^+$. Similarly, let $\Omega^-$ be the component of $\BHH-P^-$ such that $\PI \Omega^- = R^-$. Let $\X= \BHH- (\Omega^+\cup\Omega^-\cup \Delta)$. Hence, $\X$ is a mean convex domain in $\BHH$ with $\PI \X = \SI - int(R^+\cup R^- \cup Q)$. Hence, $\Gamma\subset \PI \X$.

Now, let $\mathbf{B}_n(0)$ be the ball of radius $n$ in $\BH^2$ with center $0$. Let $\mathbf{D}_n=\mathbf{B}_n(0) \times [-2m,2m]$. Let $\wh{\mathbf{D}}_n= \mathbf{D}_n\cap \X$. Let $\Gamma_n$ be the radial projection of $\Gamma$ to $\partial \wh{\mathbf{D}}_n$. Let $S_n$ be the area minimizing surface in $\wh{\mathbf{D}}_n$ with $\partial S_n = \Gamma_n$. Since $\wh{\mathbf{D}}_n$ is mean convex, $S_n$ is a smooth embedded surface in $\wh{\mathbf{D}}_n$. Again by using Lemma \ref{convergence}, we get an area minimizing surface $S$ in $\X$. By using similar ideas in Theorem \ref{APP} - Step 1, it can be showed that $\PI S = \Gamma$. While $S$ is an area minimizing surface in $\X$, it is only a minimal surface in $\BHH$. The proof follows.
\end{pf}

\begin{rmk} Notice that for smaller choice of $h_0>0$ in the Theorem above, one needs to choose the height $2m$ of the rectangles large, and the distance $s$ of the rectangles small by the construction; see Figure \ref{barrierfig}.

Recently, Kloeckner and Mazzeo has also studied these curves more extensively in \cite{KM}, where they call these curves {\em butterfly curves}. In \cite{KM}, they also constructed different families complete minimal surfaces in $\BHH$. Furthermore, they studied the asymptotic behavior of the minimal surfaces in $\BHH$.
\end{rmk}

\begin{rmk} Recently, we were able to show that when we choose $h_0<\pi$ and $2m>\pi$ sufficiently close, the butterfly curve $\Gamma_{h_0}^m$ constructed above does not bound any minimal surface, either. This example is the first non-fillable example in $\Si$ with no thin tail \cite{Co3}.  	
\end{rmk}

\section{Final Remarks} \label{secremarks}

\subsection{Infinite Curves} \

\vspace{.2cm}

In this paper, we only dealt with the finite curves, i.e. $\Gamma\subset \Si$. On the other hand, the infinite curve case is also very interesting ($\Gamma\cap (\overline{\BH^2}\times\{\pm\infty\})\neq \emptyset$). In \cite{Co2}, we studied this problem, and gave a fairly complete solution in the strongly fillable case. Kloeckner and Mazzeo studied this problem in \cite{KM}, and constructed a rich and interesting families of fillable infinite curves.

On the other hand, strong fillability question, and fillability questions are quite different in both finite and infinite curve case. While we gave a classification result for strongly fillable infinite curves in \cite{Co2}, the examples in \cite[Section 4]{Co2} shows that there are many fillable and non-fillable infinite curve families, and we are far from classification of these infinite curves in the fillable case.

\subsection{Fillable Curves} \

\vspace{.2cm}

In Section \ref{secAPPmin}, when we relax the question from "existence of area minimizing surfaces" to "existence of minimal surfaces", we see that the picture completely changes. While Theorem \ref{APP} shows that if $h(\Gamma)<\pi$, there is no area minimizing surface $\Sigma$ in $\BHH$ with $\PI \Sigma = \Gamma$,   we constructed many examples of short fillable curves $\Gamma$ in $\Si$ for any height in Section \ref{secAPPmin}.

Again, by Sa Earp and Toubiana's nonexistence theorem (Lemma \ref{nonexist}), if $\Gamma$ contains a thin tail, then there is no minimal surface $S$ in $\BHH$ with $\PI S=\Gamma$. Hence, the following classification problem is quite interesting and wide open.

\vspace{.2cm}

\noindent {\bf Classification of Fillable Curves in $\BHH$:} {\em For which curves $\Gamma$ in $\SI$, there exists a minimal surface $S$ in $\BHH$ with $\PI S=\Gamma$.}

\vspace{.2cm}

Note that Kloeckner and Mazzeo have studied this problem, and constructed many families of examples. They have also studied the asymptotic behavior of these complete minimal surfaces in $\BHH$ in \cite{KM}. Furthermore, in \cite{FMMR}, the authors have recently studied the existence of vertical minimal annuli in $\BHH$, and gave a very interesting classification. Note also that we constructed the first examples of non-fillable finite curves with no thin tail in \cite{Co3}.

\section{Appendix} \label{secappendix}

In this part, we give some technical lemmas used in the proof of the main theorem.


\subsection{Area Minimizing Catenoids in $\BHH$} \label{appendix-catenoid} \ 

\vspace{.2cm}

In this section, we study the family of minimal catenoids $\C_d$ described in \cite{NSST}, and show that for sufficiently large $d>0$,  a compact subsurface $S_d\subset\C_d$ near girth of the catenoid $\C_d$ is an area minimizing surface.

First, we recall some results on the rotationally symmetric minimal catenoids $\C_d$ \cite[Prop.5.1]{NSST}. Let $(\rho, \theta, z)$ represents the coordinates on $\BHH$ with the metric $ds^2=d\rho^2+\sinh{\rho}d\theta^2+dz^2$. Then

$$\C_d=\{(\rho,\theta,\pm\lambda_d(\rho)) \ | \ \rho\geq \sinh^{-1}d \} \  \mbox{with} \
\lambda_d(\rho)= \int_{\sinh^{-1}d}^\rho \frac{d}{\sqrt{\sinh^2{x}-d^2}}dx$$

The catenoid $\C_d$ is obtained by rotating the generating curve $\gamma_d$ about $z$-axis where $\gamma_d=\{(\rho,0,\pm\lambda_d(\rho)) \ | \ \rho\geq \sinh^{-1}d \}$. Here, $\sinh^{-1}d$ is the distance of the rotation axis to the catenoid $\C_d$, i.e. the necksize of $\C_d$.

\begin{figure}[b]

\relabelbox  {\epsfxsize=3.5in

\centerline{\epsfbox{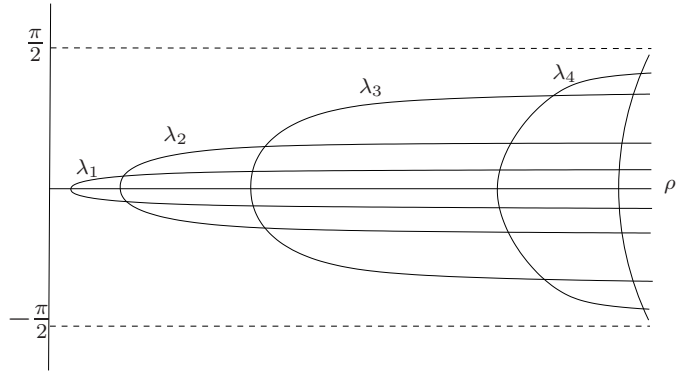}}}

\relabel{1}{$\frac{\pi}{2}$}
\relabel{2}{$-\frac{\pi}{2}$}
\relabel{3}{\tiny $\lambda_1$}
\relabel{4}{\tiny $\lambda_2$}
\relabel{5}{\tiny $\lambda_3$}
\relabel{6}{\tiny $\lambda_4$}
\relabel{7}{\tiny $\rho$}

\endrelabelbox

\caption{\label{lambda} \small For $d_i<d_{i+1}$, $\lambda_i$ represents the graphs of functions $\lambda_{d_i}(\rho)$ which are generating curves for the minimal catenoids $\C_d$. If $h(d)=\lim_{\rho\to\infty}\lambda_d(\rho)$, then $h(d)$ is monotone increasing with $h(d)\nearrow \frac{\pi}{2}$ as $d\to \infty$.}

\end{figure}

On the other hand, the asymptotic boundary of the catenoid $\C_d$ is the a pair of circles of height $\pm h(d)$, i.e. $\PI \C_d= \si\times\{-h(d),+h(d)\}\subset \Si$. Here, $h(d)=\lim_{\rho \to\infty} \lambda_d(\rho)$. By \cite{NSST}, $h(d)$ is monotone increasing function with $h(d)\searrow 0$ when $d\searrow 0$, and $h(d)\nearrow \frac{\pi}{2}$ when $d\nearrow \infty$. Hence, for any $d>0$, the catenoid $\C_d$ has height $2h(d)<\pi$ (See Figure \ref{lambda}).


By Theorem \ref{APP}, we know that the minimal catenoid $\C_d$ is not area minimizing as $\PI\C_d$ is a short curve. However, we claim that for sufficiently large $d>0$, the compact subsurfaces near the girth of $\C_d$ is indeed area minimizing. In particular, we prove the following:

\begin{lem} \label{cat1} Let $S_d^{\rho}=\C_d\cap \BH^2\times[-\lambda_d(\rho),+\lambda_d(\rho)]$ be a compact subsurface of $\C_d$. Then, for sufficiently large $d>0$, there is a $\wh{\rho}(d)>\sinh^{-1}d$ such that $S_d^{\wh{\rho}(d)}$ is an area minimizing surface.
\end{lem}

\begin{pf} Consider the upper half of the minimal catenoid $\C_d$ with the following parametrization, $\varphi_d(\rho,\theta)= (\rho, \theta,\lambda_d(\rho))$ where $\rho\geq \sinh^{-1}d$. Hence, the area of $S_d^\rho$ can be written as $$|S_d^{\rho_0}|= 2\int_0^{2\pi}\int_{\sinh^{-1}d}^{\rho_0}  \sinh{x}\sqrt{1+\frac{d^2}{\sinh^2{x}-d^2}}\ dxd\theta$$

Notice that $\partial S_d^{\rho_o}= \gamma^+_{d,\rho_o} \cup \gamma^-_{d,\rho_o}$ is a pair of round circles of radius $\rho_o$ in $\C_d$ where $\gamma^\pm_{d,\rho_o}=  \{(\rho_o,\theta, \pm\lambda_d(\rho_o)) \ | \ 0\leq \theta\leq 2\pi\}$. By \cite{NSST}, only minimal surfaces bounding $\gamma^+_{d,\rho_o} \cup \gamma^-_{d,\rho_o}$ in $\BHH$ are subsurfaces of minimal catenoids $\C_d$ and a pair of closed horizontal disks $D^+_{d,\rho_o}\cup D^-_{d,\rho_o}$ where $D^\pm_{d,\rho_o}= \{ (\rho, \theta, \pm\lambda_d(\rho_o)) \ | \ 0\leq \rho\leq \rho_o \ , \ 0\leq \theta\leq 2\pi\}$. In other words, $D^\pm_{d,\rho_o}$ is an hyperbolic disk of radius $\rho_0$ with $z=\pm\lambda_d(\rho_o)$ in $\BHH$. Recall that the area of an hyperbolic disk of radius $\rho$ is equal to $2\pi(\cosh{\rho}-1)$.

Hence, if we can show that $|S_d^\rho|<2|D_{\rho}|=4\pi(\cosh{\rho}-1)$ for some $\rho> \sinh^{-1}d$, this implies $S_d^\rho\subset \C_d$ is an area minimizing surface in $\BHH$, and we are done. Hence, we claim that there is a $\wh{\rho}(d)> \sinh^{-1}d$ such that $|S_d^\rho|<2|D_{\rho}|=4\pi(\cosh{\rho}-1)$ where $\sinh^{-1}d<\rho<\wh{\rho}(d)$. In other words, we claim the following inequality:
$$I=\int_{\sinh^{-1}d}^\rho \sinh{x}\sqrt{1+\frac{d^2}{\sinh^2{x}-d^2}}dx <\cosh{\rho} -1$$

Now, we separate the integral into two parts: $\int_{\sinh^{-1}d}^\rho \ = \int_{\sinh^{-1}d}^{\sinh^{-1}{(d+1)}}\ + \int_{\sinh^{-1}{(d+1)}}^\rho$, i.e. $I=I_1+I_2$

For the first part, clearly $$I_1=\int_{\sinh^{-1}d}^{\sinh^{-1}{(d+1)}}\sinh{x}\sqrt{1+\frac{d^2}{\sinh^2{x}-d^2}}dx <d+1\int_{\sinh^{-1}d}^{\sinh^{-1}{(d+1)}}\frac{\sinh{x}\ dx}{\sqrt{\sinh^2{x}-d^2}}$$

Recall that $\cosh^2{x}-\sinh^2{x}=1$. By substituting $u=\cosh{x}$, we get $$I_1< (d+1)\int_{\sqrt{1+d^2}}^{\sqrt{1+(d+1)^2}} \frac{du}{\sqrt{u^2-(1+d^2)}}=(d+1)\log{[u+\sqrt{u^2-(1+d^2)}]}|_{\sqrt{1+d^2}}^{\sqrt{1+(d+1)^2}}$$

This implies $$I_1<(d+1)\log{\frac{\sqrt{1+(d+1)^2}+ \sqrt{2d+1}}{\sqrt{1+d^2}}}$$.

For large $d>>0$, we obtain
$$(d+1)\log{\frac{\sqrt{1+(d+1)^2}+ \sqrt{2d+1}}{\sqrt{1+d^2}}}
\sim (d+1)\log\frac{(\sqrt{d}+\frac{1}{\sqrt{2}})^2}{d}$$

where $f(d)\sim g(d)$ represents $\frac{f(d)}{g(d)}\to 1$. After substituting $s=\sqrt{d}$ in the expression above, we get

$$\sim 2s^2 \log{\frac{s+ 1 / \sqrt{2}}{s}} \sim 2s\log(1+\frac{1}{\sqrt{2}s})^s\sim {\sqrt{2}}s={\sqrt{2d}}$$

Hence, $I_1<{\sqrt{2d}}$ for large $d>>0$.

For the second integral, we have $I_2= \int_{\sinh^{-1}{(d+1)}}^\rho \sinh{x}\sqrt{1+\frac{d^2}{\sinh^2{x}-d^2}}dx$.
Notice that the integrand $\sinh{x}\sqrt{1+\frac{d^2}{\sinh^2{x}-d^2}}=\frac{\sinh^2{x}}{\sqrt{\sinh^2{x}-d^2}}$. Hence, as $\sinh{x}<\frac{e^x}{2}$ and $\sinh^2{x}>\frac{e^{2x}-2}{4}$, we obtain $$\int\frac{\sinh^2{x}}{\sqrt{\sinh^2{x}-d^2}}dx< \int \frac{e^{2x}}{2\sqrt{e^{2x}-(2+4d^2)}}dx=\frac{\sqrt{e^{2x}-(2+4d^2)}}{2}$$.

As $\sinh^{-1}y=\log{(y+\sqrt{1+y^2})}$, after cancellations, we get $$I_2<\frac{\sqrt{e^{2\rho}-(2+4d^2)} - \sqrt{8d+2}}{2}\sim \frac{\sqrt{e^{2\rho}-(2+4d^2)}}{2}-\sqrt{2d}$$

This implies for large $d>>0$ $$I=I_1+I_2<\sqrt{2d}+(\frac{\sqrt{e^{2\rho}-(2+4d^2)}}{2}-\sqrt{2d})= \frac{\sqrt{e^{2\rho}-(2+4d^2)}}{2}$$

Now by taking $\rho=\frac{3}{2}\log{d}$ for large $d>>0$, we obtain $$I< \frac{\sqrt{e^{2\rho}-(2+4d^2)}}{2} \sim \frac{\sqrt{d^3-(2+4d^2)}}{2}\sim \frac{(d^\frac{3}{2}-2\sqrt{d}) }{2}\sim \frac{d^\frac{3}{2}}{2}- \sqrt{d}$$



On the other hand, $$\cosh{\rho}-1 = \cosh(\frac{3}{2}\log{d})-1= \frac{d^\frac{3}{2}+d^{-\frac{3}{2}}}{2}-1 \sim \frac{d^{\frac{3}{2}}}{2} $$

This shows that for $\wh{\rho}(d)=\frac{3}{2}\log{d}$, $I<\cosh{\wh{\rho}}$ and hence $|S_d^{\wh{\rho}(d)}|<2|D_{\wh{\rho}(d)}|$. Hence, the compact catenoid $S_d^{\wh{\rho}(d)}$ is an area minimizing surface in $\BHH$. The proof follows.
\end{pf}

\begin{rmk} Notice that in the lemma above, for $\wh{\rho}(d)$ is about $\frac{3}{2}$ times the neck radius of the catenoid $\C_d$, we showed that the compact slice $S_d^{\wh{\rho}(d)}$ in  $\C_d$ is an area minimizing surface. However, the comparison between $\frac{\sqrt{e^{2\rho}-(2+4d^2)}}{2}$ and $\cosh{\rho}$ indicates that if $\rho_0$ is greater than twice the neck radius of the catenoid $\C_d$ (i.e. $\rho_0>2\log(d)$), the  estimates above become very delicate, and $S_d^{\rho_0}$  is no longer area minimizing. See Remark \ref{intersectingcats} for further discussion. Note also that any subsurface of an area minimizing surface is automatically area minimizing. So, the for any $\sinh^{-1}(d)<\rho<\wh{\rho}(d)$, $S^\rho_d$ is also an area minimizing surface.
\end{rmk}

Now, we  show that as $d\to\infty$ the height $2\wh{h}(d)$ of the compact area minimizing catenoids $S_d^{\wh{\rho}(d)}$ goes to $\pi$, i.e. $\wh{h}(d)\to \frac{\pi}{2}$.

\begin{lem} \label{cat2} Let $\wh{h}(d)=\lambda_d(\wh{\rho}(d))$. Then, $\lim_{d\to\infty}\wh{h}(d)=\frac{\pi}{2}$.
\end{lem}

\begin{pf} By \cite[Prop 5.1]{NSST}, $$\lim_{d\to\infty} \wh{h}(d)=\int_0^{s(\wh{\rho}(d))}\frac{dt}{\cosh{t}}$$

By the same proposition, $s(\rho)= \cosh^{-1}(\frac{\cosh{\rho}}{\sqrt{1+d^2}})$. As $\wh{\rho}(d)=\frac{3}{2}\log{d}$, then $s(\wh{\rho}(d))\sim \sqrt{d}$. This implies $$\lim \wh{h}(d)=\int_0^\infty \frac{dt}{\cosh{t}}= \int_0^\infty \frac{du}{u^2+1}=\frac{\pi}{2}$$
\end{pf}

\begin{rmk} Notice that this lemma implies that for any height $h_o\in (0,\pi)$, there exists an area minimizing compact catenoid $S_d^\rho$ of height $h_o$. In other words, for any $h_o\in(0,\pi)$, there exists $d>0$ with $h(d)>h_o$ such that $\C_d\cap\BH^2\times [-\frac{h_o}{2},\frac{h_o}{2}]$ is an area minimizing compact catenoid in $\BHH$. Recall also that any subsurface of area minimizing surface is also area minimizing.
\end{rmk}

\begin{rmk} \label{intersectingcats} [Pairwise Intersections of Minimal Catenoids $\{\C_d\}$] \

With these results on the area minimizing subsurfaces $S_d^\rho$ in the minimal catenoids $\C_d$ in the previous part, a very interesting point deserves a brief discussion. Notice that by definition \cite{NSST}, for $d_1<d_2$, the graphs of the monotone increasing functions $\lambda_{d_1}:[\sinh^{-1}d_1,\infty)\to [0,h(d_1))$ and $\lambda_{d_2}:[\sinh^{-1}d_2,\infty)\to [0,h(d_2))$ intersect at a unique point $\rho_o\in(\sinh^{-1}d_2,\infty)$, i.e. $\lambda_{d_1}(\rho_o)=\lambda_{d_2}(\rho_o)$ (See Figure \ref{lambda}).

This implies the minimal catenoids $\C_{d_1}$ and $\C_{d_2}$ intersects at two round circles $\alpha^\pm$ of radius $\rho_o$, where $\alpha^\pm=(\rho_o,\theta,\pm\lambda_{d_1}(\rho_o))$, i.e. $\C_{d_1}\cap\C_{d_2}=\alpha^+\cup \alpha^-$.

Recall the well-known fact that two area minimizing surfaces with disjoint boundaries cannot "separate" a compact subsurface from interiors of each other. In other words, let $\Sigma_1$ and $\Sigma_2$ be two area minimizing surfaces with disjoint boundaries. If $\Sigma_1-\Sigma_2$ has a compact subsurface $S_1$ with $\partial S_1\cap\partial \Sigma_1=\emptyset$ and similarly $\Sigma_2-\Sigma_1$ has a compact subsurface $S_2$ with $\partial S_2\cap\partial \Sigma_2=\emptyset$, then $\Sigma_1'=(\Sigma_1-S_1)\cup S_2$ is an area minimizing surface with a singularity along $\partial S_1$, which contradicts to the regularity of area minimizing surfaces (Lemma \ref{AMSexist}).

This argument shows that if both $\C_{d_1}$ and $\C_{d_2}$ were area minimizing surfaces, then they must be disjoint. Hence, both $\C_{d_1}$ and $\C_{d_2}$ cannot be area minimizing surfaces at the same time. In particular, the compact area minimizing surfaces $S_{d_1}^{\rho_1}\subset \C_{d_1}$ and $S_{d_2}^{\rho_2}\subset \C_{d_2}$ must be disjoint, too.

This observation suggest an upper bound for $\wh{\rho}(d)$ we obtained in the previous part. Let $\iota(d)$ be the {\em intersection number} for $\C_d$ defined as follows: $$\iota(d)=\inf_{t>d}\{ \rho_t \ | \ \lambda_d(\rho_t)=\lambda_t(\rho_t)\}=\sup_{t<d}\{ \rho_t \ | \ \lambda_d(\rho_t)=\lambda_t(\rho_t)\}$$

The discussion above implies that $\wh{\rho}(d)<\iota(d)$ as the area minimizing surfaces $S_{d_1}^{\rho_1}\subset \C_{d_1}$ and $S_{d_2}^{\rho_2}\subset \C_{d_2}$ must be disjoint.
\end{rmk}

\subsection{Asymptotic Regularity for Tall Curves:} \

\vspace{.2cm}

The techniques in Theorem \ref{APP} also provide an elementary proof for the following result.

\begin{lem} \label{surface} Let $\Sigma$ be a complete area minimizing surface in $\BHH$. Let $\overline{\Sigma}$ be the closure of $\Sigma$ in $\overline{\BHH}$, and let $\Gamma=\PI\Sigma$. If $\Gamma$ is a tall curve, then $\overline{\Sigma}$ is a surface with boundary.
\end{lem}

\begin{pf} By Step 1 in Theorem \ref{APP}, and Corollary \ref{MCHcor},  $\Sigma\subset MCH(\Gamma)$. This will imply $\overline{\Sigma}=\Sigma\cup \Gamma$ as follows. Let $p$ be point in $\Gamma$, and let $U_p$ be a sufficiently small neighborhood of $p$ in $\overline{\BHH}$ so that $U_p\cap \Gamma$ is a small arc $\gamma_p$ in $\Gamma$. We claim that $\overline{\Sigma}\cap U_p$ is an embedded surface with boundary.

There are two cases. Either $\gamma_p$ is vertical, or not. Assume that $\gamma_p$ is not a vertical segment. Consider the upper half space model for $\BHH$, and without loss of generality, let $p=(0,0,0)$. As in Step 2 of Theorem \ref{APP}, let $\wh{\psi}_t(x,y,z)= (tx,ty,z)$ be the isometry of $\BHH$. Let $\Sigma_n=\wh{\psi}_n(\Sigma)$. Then, the sequence $\{\Sigma_n\}$ has a subsequence converging to $\wh{\Sigma}$ with $\PI \wh{\Sigma}=\wh{\Gamma}$ by Lemma \ref{convergence}. Let $\wh{\gamma}$ be the component of $\wh{\Gamma}$ containing $p$. Since $\gamma_p$ is not vertical, $\wh{\gamma}$ is either the straight line $l=\BR\times\{0\}\times\{0\}$ ($x$-axis) or a half line $l^\pm=\BR^\pm\times\{0\}\times\{0\}$. Since $\Gamma$ is tall, $\Gamma$ does not have any thin tail (Definition \ref{tail}). This excludes the case $\wh{\gamma}$ is a half line, and hence $\wh{\gamma}=l$. In particular, in cylindrical model, this shows $\wh{\gamma}=S^1_\infty\times \{0\}$. As $\Gamma$ is tall, by the proof of \cite[Lemma 8.6]{Co1} (the last paragraph), the component $T$ of $\wh{\Sigma}$ with $\PI T\supset \wh{\gamma}$ is the horizontal plane $\BH^2\times\{0\}$, i.e. $T=\BH^2\times\{0\}$. Note that the isometry $\wh{\psi}_n$ fixes the $z$-direction. As $\Sigma_n\to\wh{\Sigma}$, this proves that $\Sigma_n$ is graph over $\BH^2\times\{0\}$ near $p$ for sufficiently large $n$. As $\Sigma_n$ is isometric image of $\Sigma$, this shows that $\overline{\Sigma}\cap U_p$ is an embedded surface with boundary.

In the second case, we  assume $\gamma_p$ is vertical. Again, assume $\alpha_p=\{0\}\times\{0\}\times [-c,+c]$. Similarly, define $\Sigma_n=\wh{\psi}_n(\Sigma)$, $\wh{\Sigma}$, and $\wh{\Gamma}$ as before. Let $l_c^+=\BR^+\times\{0\}\times\{c\}$, and $l_c^-=\BR^-\times\{0\}\times\{c\}$. Similarly define $l_{-c}^\pm$. By construction,  $\wh{\Gamma}$ is either $l_c^+\cup l_{-c}^+$ or $l_c^-\cup l_{-c}^-$, say Case A, or $l_c^+\cup l_{-c}^-$ or $l_c^-\cup l_{-c}^+$, say Case B. The first two cases (Case A) are similar, and last two cases (Case B) are similar.

In cylindrical model, $\wh{\Gamma}$ correspond to the tall rectangle $R=[0,\pi]\times[-c,c]$ (or $R=[\pi,2\pi]\times[-c,c]$)  in Case A. In this case, as $\Gamma$ is tall, this shows that $\gamma_p$ belongs to a vertical segment $\alpha_p$ longer than $\pi$ in $\Gamma$, i.e. $c>\frac{\pi}{2}$. There is a unique area minimizing surface $P$ with $\PI P=R$ by Lemma \ref{rectangle}. This proves that $\wh{\Sigma}=P$ and hence $\overline{\Sigma}\cap U_p$ is an embedded surface with boundary as before.
In Case B, $\wh{\Gamma}$ correspond to the curve $R'=(\{0\}\times[-c,c])\cup ([0,\pi]\times\{c\})\cup(\{\pi\}\times[-c,c])\cup ([\pi,2\pi]\times\{-c\}$. Again, $R'$ bounds a unique area minimizing surface, which is a graph over $\BH^2\times\{0\}$ by \cite[Proposition 2.1 (3)]{ST}. The proof follows.
\end{pf}

\begin{rmk} Note that Kloneckner and Mazzeo proved higher order ($\C^{k,\alpha}$) asymptotic regularity for embedded minimal surfaces in $\BHH$ by using different techniques in \cite[Section 3]{KM}.
\end{rmk}


\end{document}